\documentclass[12pt]{article}
\usepackage{amssymb,latexsym,amsmath}

\begin{document}
\pagestyle{plain} \headheight=5mm  \topmargin=-5mm

\title{Infinitely generated Lawson homology groups on some rational
projective varieties}
\author{ Wenchuan Hu }
%\author{Huitao Feng , Wenchuan Hu \quad and\quad Weiping Zhang\\
%Department of Mathematics,\\
%SUNY at Stony Brook, Stony Brook, NY, 11794-3651 }
%\date{April 24, 2004}

\maketitle
\newtheorem{Def}{Definition}[section]
\newtheorem{Th}{Theorem}[section]
\newtheorem{Prop}{Proposition}[section]
\newtheorem{Not}{Notation}[section]
\newtheorem{Lemma}{Lemma}[section]
\newtheorem{Rem}{Remark}[section]
\newtheorem{Cor}{Corollary}[section]

\def\s{\section}
\def\ss{\subsection}

\def\d{\begin{Def}}
\def\t{\begin{Th}}
\def\p{\begin{Prop}}
\def\n{\begin{Not}}
\def\la{\begin{Lemma}}
\def\r{\begin{Rem}}
\def\c{\begin{Cor}}
\def\ee{\begin{equation}}
\def\aa{\begin{eqnarray}}
\def\ya{\begin{eqnarray*}}
\def\bd{\begin{description}}

\def\ed{\end{Def}}
\def\et{\end{Th}}
\def\epo{\end{Prop}}
\def\en{\end{Not}}
\def\el{\end{Lemma}}
\def\er{\end{Rem}}
\def\ec{\end{Cor}}
\def\eee{\end{equation}}
\def\eaa{\end{eqnarray}}
\def\ey{\end{eqnarray*}}
\def\ebd{\end{description}}

\def\nn{\nonumber}
\def\bp{{\bf Proof.}\hspace{2mm}}
\def\qe{\hfill$\Box$}
\def\lj{\langle}
\def\rj{\rangle}
\def\dd{\diamond}
\def\ox{\mbox{}}
\def\lb{\label}
\def\rel{\;{\rm rel.}\;}
\def\vp{\varepsilon}
\def\ep{\epsilon}
\def\mod{\;{\rm mod}\;}
\def\exp{{\rm exp}\;}
\def\Lie{{\rm Lie}}
\def\dim{{\rm dim}}
\def\im{{\rm im}\;}
\def\Lag{{\rm Lag}}
\def\Gr{{\rm Gr}}
\def\span{{\rm span}}
\def\Spin{{\rm Spin}}
\def\sign{{\rm sign}\;}
\def\Supp{{\rm Supp}\;}
\def\Sp{{\rm Sp}\;}
\def\ind{{\rm ind}\;}
\def\rank{{\rm rank}\;}
\def\Sg{{\Sp(2n,\C)}}
\def\Na{{\cal N}}
\def\det{{\rm det}\;}
\def\dist{{\rm dist}}
\def\deg{{\rm deg}}
\def\Tr{{\rm Tr}\;}
\def\ker{{\rm ker}\;}
\def\Vect{{\rm Vect}}
\def\H{{\bf H}}
\def\K{{\rm K}}
\def\R{{\bf R}}
\def\C{{\bf C}}
\def\Z{{\bf Z}}
\def\N{{\bf N}}
\def\F{{\bf F}}
\def\Da{{\bf D}}
\def\A{{\bf A}}
\def\La{{\bf L}}
\def\x{{\bf x}}
\def\y{{\bf y}}
\def\Ga{{\cal G}}
\def\Ha{{\cal H}}
\def\L{{\cal L}}
\def\Pa{{\cal P}}
\def\Ua{{\cal U}}
\def\E{{\rm E}}
\def\J{{\cal J}}

\def\m{{\rm m}}
\def\ch{{\rm ch}}
\def\gl{{\rm gl}}
\def\Gl{{\rm Gl}}
\def\Sp{{\rm Sp}}
\def\sf{{\rm sf}}
\def\U{{\rm U}}
\def\O{{\rm O}}
\def\F{{\rm F}}
\def\P{{\rm P}}
\def\D{{\rm D}}
\def\T{{\rm T}}
\def\Sa{{\rm S}}

%\begin{center}{\bf Table of Contents}\end {center}

\begin{center}{\bf \tableofcontents}\end {center}

\begin{center}{\bf Abstract}\end {center}
{\hskip .2 in} We construct rational projective 4-dimensional
varieties with the property that certain Lawson homology groups
tensored with ${ \mathbb{Q}}$ are infinite dimensional
${\mathbb{Q}}$-vector spaces. More generally, each pair of integers
$p$ and $k$, with $k\geq 0$, $p>0$, we find a projective variety
$Y$, such that $L_pH_{2p+k}(Y)$ is infinitely generated.

We also construct two singular rational projective 3-dimensional
varieties $Y$ and $Y'$ with the same homeomorphism type but
different Lawson homology groups, specifically $L_1H_3(Y)$ is not
isomorphic to $L_1H_3(Y')$ even up to torsion.

%Based on a result of Clemens that the Griffiths groups (a special
%case of Lawson homology) of 1-cycles on a certain 3-fold, tensored
%with ${\mathbb{Q}}$, is an infinite dimensional ${\mathbb{Q}}$-vector space,

\s {Introduction} {\hskip .2 in} This paper gives examples of
singular rational projective 4-dimensional varieties with infinitely
generated Lawson homology groups even modulo torsion. This is
totally different from the smooth case ([Pe], also [H1]), where it
is known that all Lawson homology groups of rational fourfolds are
finitely generated.

This paper also gives examples of singular rational projective
3-dimensional varieties with the same homeomorphism type but
different Lawson homology groups.

\medskip
For an algebraic variety $X$ over $\mathbb{C}$, the {\bf Lawson
homology} $L_pH_k(X)$ of $p$-cycles is defined by

$$L_pH_k(X) := \pi_{k-2p}({\cal Z}_p(X)), \quad k\geq 2p\geq 0$$
where ${\cal Z}_p(X)$ is provided with  a natural topology. For
general background, the reader is referred to the survey paper
[L2].

Clemens showed that the {\bf Griffiths group} of 1-cycles (which is
defined to be  the group of algebraic 1-cycles homologically
equivalent to zero modulo l-cycles algebraically equivalent to zero)
may be infinitely generated even modulo the torsion elements for
general quintic hypersurfaces in $\P^4$ (cf.[C]).  Friedlander
showed that $L_1H_2(X)$ is exactly the algebraic 1-cycles modulo
algebraic equivalence (cf. [F]). Hence the Griffiths group of
1-cycles for $X$ is a subgroup of $L_1H_2(X)$.

\medskip
This leads to the following question:

\begin{enumerate}
\item[\bf (Q):] Can one show that $L_pH_{2p+j}(X)$ is {\bf not}
finitely generated for some projective variety $X$ where $j>0$ ?
\end{enumerate}

In this paper we shall construct, for any given integers $p$ and
$j>0$, examples of rational varieties $X$ for which
$L_pH_{2p+j}(X)$, as an abelian group, is infinitely generated.
Thus, we answer affirmatively the above question :

 {\Th There exists  rational projective variety $X$ with $\dim(X)=4$ such that
$L_1H_3(X)\otimes {\mathbb{Q}}$ is {\bf not} a finite dimensional
${\mathbb{Q}}$-vector space. }

%By the projective bundle theorem given by Friedlander and Gabber([FG]),
% and the concrete computation of the resolution of
%singularities (the singularity we will meet is of the so-called
%$A_n$ type), we obtain the following:

%{\Th(?) There exists  smooth rational projective 4-manifold $X$
%such that $L_1H_3(X)\otimes$ is {\bf NOT} a finite dimensional
%${\mathbb{Q}}$-vector space. }

\medskip
By using the projective bundle theorem given by Friedlander and
Gabber([FG]), we have the following corollary:

{\Cor For any $p\geq 1$, there exists projective algebraic variety
$X$ such that $L_pH_{2p+1}(X)$ is {\bf not} a finitely generated
abelian group.}

\medskip
More generally, we have

 {\Th For integers $p$ and $k$, with $k\geq
0, p> 0$, we can find a projective variety $Y$,  such that
$L_pH_{2p+k}(Y)$ is infinitely generated.}

{\Rem The smoothness is essential here. Compare Theorem 1.1 with the
following result proved by C. Peters.}

{\Th ([Pe]) For any smooth projective variety $X$ over $\mathbb{C}$
with ${\rm Ch}_0(X)\otimes{\mathbb{Q}}\cong {\mathbb{Q}}$, the
natural map $\Phi: L_1H_*(X)\otimes{\mathbb{Q}}\rightarrow
H_*(X,{\mathbb{Q}})$ is injective. In particular,
$L_1H_*(X)\otimes{\mathbb{Q}}$ is a finite dimensional
${\mathbb{Q}}$-vector space.}

\medskip
Any rational variety $X$ (smooth or not) has the property that ${\rm
Ch}_0(X)\otimes {\mathbb{Q}}\cong {\mathbb{Q}}$.

\medskip
Applying the same construction to hypersurfaces in $\P^3$, we obtain
the following:

 {\Th There exist two rational 3-dimensional projective varieties
 $Y$ and $Y'$ which are homeomorphic  but for which the Lawson homology
 groups $L_1H_3(Y,{\mathbb{Q}})$ and $L_1H_3(Y',{\mathbb{Q}})$ are not isomorphic even
up to torsion. }

{\Rem In fact, these varieties in Theorem 1.4 have exactly
\emph{one} isolated singular point.}

% I would like to express my gratitude to my advisor, Blaine
%Lawson, for proposing this problem and for all his help.

%section 222222222222222222222222222**************************************************************************
 \s{Lawson Homology }
{\hskip .2 in} In this section we briefly review the definitions and
results used in the next section.
 Let $X$ be a projective variety of dimension $m$ over $\mathbb{C}$. The
 group of $p$-cycles on $X$ is the free abelian group ${\cal
 Z}_p(X)$ generated by irreducible $p$-dimensional subvarieties.

{\Def The \textbf{Lawson homology} $L_pH_k(X)$ of $p$-cycles on
$X$ is defined by
$$L_pH_k(X) :=\pi_{k-2p}({\cal Z}_p(X)), k\geq 2p\geq 0,$$
where ${\cal Z}_p(X)$ is provided with a natural, compactly
generated topology (cf. [F], [L1], [L2]).}

{\Def  The {\bf Griffiths group} ${\rm Griff}_p(X)$ of $p$-cycles
on $X$ is defined by
$${\rm Griff}_p(X) :={\cal Z}_p(X)_{hom}/{\cal Z}_p(X)_{alg} $$
where ${\cal Z}_p(X)_{hom}$ denotes algebraic $p$-cycles homologous
to zero and ${\cal Z}_p(X)_{alg}$ denotes algebraic $p$-cycles which
are algebraically equivalent to zero.}

{\Rem It was shown by Friedlander that $L_pH_{2p}(X)\cong {\cal
Z}_p(X)/{\cal Z}_p(X)_{alg}$ (cf. [F]). Hence the {Griffiths
group} ${\rm Griff}_p(X)$ is a subgroup of the Lawson homology
$L_pH_{2p}(X)$. Therefore, for any projective variety $X$ (its
homology groups are finitely generated), ${\rm Griff}_p(X)$ is
infinitely generated if and only if $L_pH_{2p}(X)$ is.}

{\Rem For a quasi-projective variety $U$, $L_pH_k(U)$ is also
well-defined and independent of the projective embedding (cf. [Li],
[L2]).}

\medskip
Let $V\subset U$ be a Zariski open subset of a quasi-projective
variety $U$. Set $Z=U-V$. Then we have

{\Th([Li]) There is a long exact sequence for the pair $(U, Z)$,
i.e.,
\begin{eqnarray}
\cdots\rightarrow L_pH_k(Z)\rightarrow L_pH_k(U)\rightarrow
L_pH_k(V)\rightarrow L_pH_{k-1}(Z)\rightarrow\cdots \label{eq1}
\end{eqnarray}}

{\Rem For any quasi-projective variety $U$, $L_0H_k(U)\cong
H^{BM}_k(U)$, where $H^{BM}_k(U)$ is the Borel-Moore homology. This
follows from the Dold-Thom Theorem [DT].}

\medskip
As a direct application of this long exact sequence, one has the
following results [Li]:
\begin{enumerate}
\item[2.1] Let $U=\P^{n+1}$ and $V=\P^{n+1}-\P^n$. By the Complex
Suspension Theorem [L1], we have, $L_pH_{2n}(\mathbb{C}^n)={\bf Z}$;
$L_pH_k(\mathbb{C}^n)=0$ for any $k\neq 2n$ and $k\geq 2p\geq 0$.

\item[2.2] Let $U=\mathbb{C}^{n}$ and $V_{n-1}\subset \mathbb{C}^{n}$ be a closed
algebraic set. Set $V_n=\mathbb{C}^{n}-V_{n-1}$. Then we have
$$0 \rightarrow L_pH_{2n+1}(V_n)\rightarrow
L_pH_{2n}(V_{n-1})\rightarrow L_pH_{2n}({\mathbb{C}^n}) \rightarrow
L_pH_{2n}(V_n)\rightarrow L_pH_{2n-1}(V_{n-1})\rightarrow 0$$ and
$$L_pH_{k+1}(V_n)\cong L_pH_{k}(V_{n-1}),\quad \rm k\neq 2n,2n+1.$$
\end{enumerate}

 \s{An Elementary Construction}
{\hskip .2 in} $\bf Construction$  Let $X=(f(x_0,\cdots,x_{n+1})=0)$
be a general hypersurface in $\P^{n+1}$ with degree $d$, and let
$V_{n}:=X-X\cap\{\P^{n}=(x_0=0)\}$ be the affine part, i.e.
$V_{n}\subset \mathbb{C}^{n+1}$. Define
$V_{n+1}:={\mathbb{C}^{n+1}}-V_{n}$, then $V_{n+1}$ can be viewed as
an affine variety in ${\mathbb{C}^{n+2}}$ defined by $x_{n+2}\cdot
f(1,x_1,\cdots,x_{n+1})-1=0$, where
$V_{n}=(f(1,x_1,\cdots,x_{n+1})=0)$. Denoted by $\overline{V_{n+1}}$
the projective closure of $V_{n+1}$ in $\P^{n+2}$ and set $Z_{n} =
\overline{V_{n+1}}-V_{n+1}$.

\medskip
We leave the study of this case where $n=1$, and $X$ is a smooth
plane curve, as an exercise.

\ss {Application to the Case $n=2$} {\hskip .2 in}  In this
subsection, I will show that there exist two rational projective
3-dimensional varieties with the same singular homology groups but
different Lawson homology.

The following result  proved by Friedlander will be used several
times: {\Th (Friedlander [F1]) Let $X$ be any smooth projective
variety of dimension $n$. Then we have the following isomorphisms
$$
\left\{
\begin{array}{l}
 L_{n-1}H_{2n}(X)\cong \Z,\\
 L_{n-1}H_{2n-1}(X)\cong H_{2n-1}(X,\Z),\\
 L_{n-1}H_{2n-2}(X)\cong H_{n-1,n-1}(X,\Z)=NS(X)\\
L_{n-1}H_{k}(X)=0 \quad for\quad k> 2n.\\

\end{array}
\right.
$$}

For a finitely generated abelian group $G$, we denote by ${\rm
rk}(G)$ the rank of $G$.

\medskip
Let $X\subset \P^3$ be a general surface with degree $d=4$. Then
$V_2=X-X\cap {\P^2}$ and $C:=X\cap {\P^3}$ is a smooth curve in
${\P^2}$.

{\Lemma  ${\rm rk}(L_1H_{2}(X))\cong {\rm rk}(L_1H_{2}(V_2))+1$;
${\rm rk}(L_1H_{3}(V_2))=0$.}

\medskip
\bp Applying Theorem 2.1 to the pair $(X, C)$ and Theorem 3.1 for
$X$, we get
$$0\rightarrow L_1H_3(V_2)\rightarrow L_1H_2(C)\rightarrow
L_1H_2(X)\rightarrow L_1H_{2}(V_2)\rightarrow 0.$$

Note that $ L_1H_2(C)\cong\Z$ and the map $L_1H_2(C)\rightarrow
L_1H_2(X)$ is injective, and so we get $L_1H_{3}(V_2)=0$. Therefore,
by the above long exact sequence, we have ${\rm
rk}(L_1H_{2}(X))\cong {\rm rk}(L_1H_{2}(V_2))+1$.\qe

{\Lemma ${\rm rk}(L_1H_2(Z_2))=1$; ${\rm rk}(L_1H_3(Z_2))=6$;
${\rm rk}(L_1H_4(Z_2))= 2$.}

\medskip
\bp Note that $Z_{2} = \overline{V_{3}}-V_{3}$ is defined by
$(x_4\cdot f(0,x_1,...,x_3)=0,x_0=0)$ in $\P^4$. Let
$C'=(x_4=0)\cap(f(0,x_1,\cdots,x_3)=0)$ in the hyperplane
$(x_0=0)\subset {\P^3}$. It is easy to see that $C'\cong C$. Then
$Z_2={\P^2}\cup \Sigma_p(C)$, where $\Sigma_p(C)$ means the joint of
$C$ and the point $p=[1:0:\cdots:0]$. By applying Theorem 2.1 to the
pair $(Z_2, \Sigma C)$, we get
$$
 \cdots\rightarrow L_1H_3(Z_2-\Sigma C)\rightarrow L_1H_2(\Sigma
 C)\rightarrow L_1H_2(Z_2)\rightarrow L_1H_{2}(Z_2-\Sigma
 C)\rightarrow 0.
$$

Note that $Z_2-\Sigma C\cong \P^2-C$ and $L_1H_3(\P^2-C)=0$.
Therefore ${\rm rk}(L_1H_2(Z_2))=1$. Moreover, since
$L_1H_4(\P^2-C)=\Z$ and $L_1H_4(\Sigma C)=\Z$ and $\P^2\cap \Sigma
C=C$ is a curve. The last statement follows. Recall that the
Complex Suspension Theorem and Dold-Thom Theorem, we have
$L_1H_3(\Sigma C)\cong L_0H_1(C)\cong H_1(C)$. By assumption, $C$
is a plane curve of degree 4. The adjunction formula gives ${\rm
rk}(H_1(C))=6$. The second statement follows.
\qe

{\Lemma ${\rm rk}(L_1H_2(\overline{V_3}))\leq 1$; ${\rm
rk}(L_1H_3(\overline{V_3}))={\rm rk}(L_1H_{2}({X}))+{\rm
rk}(L_1H_2(\overline{V_3}))+4$.}

\medskip
\bp Applying Theorem 2.1 to the pair $(\overline{V_3},Z_2)$ with
$p=1$, we have
$$
0 \rightarrow L_1H_{3}({Z_2})\rightarrow
 L_1H_{3}(\overline{V_3})\rightarrow L_1H_{3}(V_{3})\rightarrow
 L_1H_{2}({Z_2})\rightarrow L_1H_{2}(\overline{V_3}) \rightarrow 0
$$
since Lemma 3.1 gives $L_1H_{2}(V_{3})=0$ and
$L_1H_{4}(V_{3})\cong L_1H_{3}(V_{2})=0$. Hence ${\rm
rk}(L_1H_{2}({\overline{V_3}}))\leq 1$. Moreover, we have ${\rm
rk}(L_1H_{3}({Z_2}))-{\rm rk}(L_1H_{3}(\overline{V_3}))+ {\rm
rk}(L_1H_{3}(V_{3}))-{\rm rk}(L_1H_{2}({Z_2}))+ {\rm
rk}(L_1H_{2}(\overline{V_3}))=0$. By Lemma 3.2, we get
$$ 6-{\rm rk}(L_1H_{3}(\overline{V_3}))+({\rm
rk}(L_1H_{2}({X}))-1)-1+{\rm rk}(L_1H_{2}(\overline{V_3}))=0
$$
\qe

{\Lemma ${\rm Sing}(\overline{V_3})\cong \{X\cap
({x_0=0})\}\cup\{p\}\cong C\cup\{p\}$.}

\medskip
\bp It follows from a direct computation. By definition,
$$
{\rm Sing}(\overline{V_3})=\left\{
\begin{array}{l}
 F(x_0,x_1,x_2,x_3,x_4)=0,\\
 dF(x_0,x_1,x_2,x_3,x_4)=0\\

\end{array}
\right\}
$$

$$
=\left\{
\begin{array}{l}
x_0^{d+1}-x_4\cdot f(x_0,x_1,x_2,x_3)=0,\\
 (d+1)x_0-x_4\cdot \frac{\partial f}{\partial x_0}=0,\\
-x_4\cdot \frac{\partial f}{\partial x_1}=0,\\
-x_4\cdot \frac{\partial f}{\partial x_2}=0,\\
-x_4\cdot \frac{\partial f}{\partial x_3}=0,\\
f(x_0,x_1,x_2,x_3)=0
\end{array}
\right\}
$$

$$
=\left\{
\begin{array}{l}
x_0=0,\\
 x_4\cdot \frac{\partial f}{\partial x_0}=0,\\
x_4\cdot \frac{\partial f}{\partial x_1}=0,\\
x_4\cdot \frac{\partial f}{\partial x_2}=0,\\
x_4\cdot \frac{\partial f}{\partial x_3}=0,\\
f(x_0,x_1,x_2,x_3)=0
\end{array}
\right\}
$$

$$=\bigg\{x_0=x_3=f(x_0,x_1,x_2)=0\bigg\}\cup \bigg\{x_0=\frac{\partial f}{\partial x_0}=\frac{\partial
f}{\partial x_1}=\frac{\partial f}{\partial x_2}=\frac{\partial
f}{\partial x_3}=0\bigg\}$$
$$\cong \{x_0=f(x_0,x_1,x_2)=0\}\cup \{p \}\cong C\cup \{p \}$$
since $C=(f=0)$ is smooth by our assumption.
 \qe\\

{\Rem Note that $p$ is an isolated singular point and  the
singularity $C=X\cap ({x_0=0})$ is of $A_n$-type. We can resolve the
singularity of this part by  blowing up twice over the singularity,
i.e., by blowing up over the singularity for the first time  and
then blowing up the singularity of the proper transform of the first
blowup. We denote by $\widetilde{\overline{V_3}}$ the proper
transform of $\overline{V_3}$ with the exceptional divisor $D_1$ for
the first blowup and $\widetilde{\widetilde{\overline{V_3}}}$ the
proper transform of $\widetilde{\overline{V_3}}$ with the
exceptional divisor $D_2$ for the second blowup. Both $D_1$ and
$D_2$ are isomorphic to a fiber bundle over $C$ with fibre the union
of two $\P^1$ intersecting at exactly one point. See the appendix
for the computation of a concrete example. }

Now $\widetilde{\widetilde{\overline{V_3}}}$ has only one singular
point, denote by $q$.

{\Lemma The singular point $q$ in
$\widetilde{\widetilde{\overline{V_3}}}$ can be resolved by one blow
up whose exceptional divisor is isomorphic to $X$.}

\medskip
\bp It follows from a trivial computation.\qe

\medskip

We denote by $W_3$ the proper transform of the blow up in the above
lemma. Note that $W_3$ is a \textbf{smooth} rational threefold. We
have the following property on $W_3$:

{\Prop For a smooth surface $X\subset \P^3$, the $W_3$ thus
constructed is a smooth rational threefold with a fixed homeomorphic
type, i.e., for two smooth surfaces $X$ and $X'$ in $\P^3$, the
corresponding smooth rational threefolds $W_3$ and $W_3'$ are
homeomorphic.}

\medskip
\bp Note that $V_3$ is a hypersurface in $\P^4$. Let
$(f_t(x_0,\cdots,x_4)=0)\subset \P^4$ be a family of hypersurface
such that $V_3=(f_0=0)$ is transversal to the hypersurface
$H=(x_0=0)$. Let $\Delta$ be a neighborhood of $t=0$ such that
$(f_t=0)$ is transversal to $H$ for all $t\in\Delta$. Let $W\subset
\P^4\times \Delta$ be the (analytic) variety defined by
$F(x,t):=f_t(x)=0$. Then we have the following incidence
correspondence

$$
\begin{array}{cccc}
 W & \subset & \P^4\times \Delta& \\
 \downarrow \pi &   & \downarrow \pi&   \\
\Delta& = & \Delta,&
\end{array}
$$
with $\pi^{-1}(t)\cap W=W^t$.

By Remark 3.1 and Lemma 3.5, we get a smooth variety $\widetilde{W}$
by blowing up twice along 2-dimensional singularity of $Sing(W)$ and
once for the remainder 1-dimensional singularity of $Sing(W)$.
Denote by $E$ the exceptional divisor of the last step. We claim
that the map $\tilde{\pi}:\widetilde{W}\rightarrow \Delta$ is a
smooth proper submersion. In fact, let $v$ be a vector field of
$\Delta$ and let $\tilde{v}$ be a lifting in
$\Gamma(\widetilde{W},T\widetilde{W})$ such that
$\tilde{\pi}_*(\tilde{v})=v$.

Denote by $\varphi_t$ (resp. $\tilde{\varphi}_t$) the flow
determined by $v$ (resp. $\tilde{v}$). Then
$\tilde{\varphi}_t:\widetilde{W}^0\rightarrow \widetilde{W}^t$ gives
the homeomorphism between two fiber of $\tilde{\pi}$ from
Ehresmann's Theorem [V]. This implies the result of the proposition.

 \qe

From this proposition, we have the following

{\Cor For all smooth surfaces $X\subset \P^3$ of fixed degree, the
$\widetilde{\widetilde{\overline{V_3}}}$ thus constructed has  a
fixed homeomorphism type.}

\medskip
\bp In this proof of the proposition, we actually can choose
$\tilde{v}$ such that 1)$\tilde{v}$ is tangent to W; 2) $\tilde{v}$
is tangent to the exceptional divisor $E$. Then the flow of
$\tilde{v}$ gives the homeomorphism of any two fibers.

\qe

We want to show that some Lawson homology group of
$\widetilde{\widetilde{\overline{V_3}}}$ may vary when the general
$X$ varies in $\P^3$.

{\Th There exist two rational 3-dimensional projective varieties
$Y$, $Y'$ such that $Y$ is homeomorphic to $Y'$ but the Lawson
homology group $L_1H_3(Y)$ is not isomorphic to $L_1H_3(Y')$ even up
to torsion. }

\medskip
\bp If $X\subset \P^3$ is a general smooth quartic surface, then the
Picard group ${\rm Pic}(X)\cong \Z$ by Noether-Lefschetz Theorem.
For details, see e.g. Voisin [V]. But it is well known that there
are still many special smooth quartic surfaces $X'$ in $\P^3$ with
${\rm rk}({\rm Pic}(X'))$ as big as 20. Note that by Theorem 3.1 and
the Weak Lefschetz Theorem $L_1H_2(X)\cong {\rm Pic}(X)$ for any
smooth surface $X$ in $\P^3$.

Now we choose smooth $X$ with $L_1H_2(X)\cong \Z$ and $X'$ with
$L_1H_2(X')\cong \Z^{20}$. Set
$Y:=\widetilde{\widetilde{\overline{V_3}}}$ and
$Y':=\widetilde{\widetilde{\overline{V_3'}}}$. Let $W_3$ (resp.
$W_3'$) be as in Proposition 3.1. From the proof of Lemma 2.1 in
[H1], we have the commutative diagram

$$
\begin{array}{ccccccccccc}
\cdots\rightarrow & L_1H_3(E) & {\rightarrow} & L_1H_3(W_3)
 & {\rightarrow} & L_1H_3(W_3-E) & {\rightarrow} & L_1H_{2}(E) & \rightarrow & \cdots & \\

 & \downarrow &   & \downarrow &   & \downarrow  \cong&   & \downarrow &   &  &\\

 \cdots\rightarrow & L_1H_3(q) & {\rightarrow} &
L_1H_3(Y) & {\rightarrow}& L_1H_3(Y-q)&{\rightarrow} & L_1H_{2}(q) &
\rightarrow & \cdots&

\end{array}
$$

 By Lemma 3.5, we know $E\cong X$. By Theorem 3.1, we  have $L_1H_3(X)\cong
 H_3(X)$. By the Lefschetz Hyperplane Theorem, we know $X$ is simply
 connected. Since $q$ is a point, we have $L_1H_3(Y)\cong L_1H_3(Y-q)\cong
 L_1H_3(W_3-E)$ and $L_1H_2(Y)\cong L_1H_2(Y-q)\cong
 L_1H_2(W_3-E)$.

The top row of the above commutative diagram turns into the long
exact sequence
$$0  {\rightarrow}  L_1H_3(W_3)
  {\rightarrow}  L_1H_3(Y)  {\rightarrow}  L_1H_{2}(X)  \rightarrow
  L_1H_2(W_3)\rightarrow L_1H_2(Y)\rightarrow 0
$$

Therefore, we have
$$ {\rm rk} L_1H_3(W_3)-{\rm rk} L_1H_3(Y)+{\rm rk} L_1H_2(X)-{\rm rk}
L_1H_2(W_3)+{\rm rk} L_1H_2(Y)=0
$$

Since $W_3$ is a smooth rational threefold, we have
$L_1H_3(W_3)\cong H_3(W_3)$, $L_1H_2(W_3)\cong H_2(W_3)$ ([FHW,
Prop. 6.16]) and by Proposition 3.1 $H_i(W_3)\cong H_i(W_3')$ for
all $i$.

$$
\begin{array}{ccc}
(*)\quad {\rm rk} L_1H_3(Y)&=& {\rm rk} H_3(W_3)+{\rm rk}
L_1H_2(X)-{\rm rk} H_2(W_3)+{\rm rk} L_1H_2(Y)\\
\end{array}
$$

By applying Theorem 2.1 to $(\widetilde{\overline{V_3}},D_1)$, we
get
$$  \cdots \rightarrow L_1H_2(D_1)\rightarrow L_1H_2(\widetilde{\overline{V_3}})\rightarrow
L_1H_2(\widetilde{\overline{V_3}}-D_1)\rightarrow 0
$$

Hence
$$
\begin{array}{ccc}
 {\rm rk}L_1H_2(\widetilde{\overline{V_3}})&\leq &{\rm
rk}L_1H_2(D_1)+{\rm rk}L_1H_2(\widetilde{\overline{V_3}}-D_1)\\
&=& {\rm rk}L_1H_2(D_1)+{\rm rk}L_1H_2(\overline{V_3}-q)\\
&=& {\rm rk}L_1H_2(D_1)+{\rm rk}L_1H_2(\overline{V_3})\\
&\leq& {\rm rk}L_1H_2(D_1)+1 \quad({\rm Lemma \quad 3.3})
\end{array}
$$

Similarly, $$
\begin{array}{ccc}
 {\rm rk}L_1H_2(\widetilde{\widetilde{\overline{V_3}}})&\leq&
 {\rm rk}L_1H_2(D_2)+{\rm rk}L_1H_2(\widetilde{\overline{V_3}})
\end{array}
$$

Therefore, $$
\begin{array}{ccc}
 {\rm rk}L_1H_2(\widetilde{\widetilde{\overline{V_3}}})&\leq&
 {\rm rk}L_1H_2(D_2)+{\rm rk}L_1H_2(D_1)+1.
\end{array}
$$

Since $D_1$ (also $D_2$) is isomorphic to a $\P^1\cup\P^1$-bundle
over a smooth curve $C$, it is easy to compute, by using Theorem 2.1
and the Projective Bundle Theorem [FG], that
$$ {\rm rk}L_1H_2(D_1)\leq {\rm rk}L_1H_2(C)+2\cdot {\rm
rk}L_0H_0(C)=1+2\times 1=3.
$$

Therefore
$${\rm
rk}L_1H_2(\widetilde{\widetilde{\overline{V_3}}})\leq 3+ 3+1=7.
$$

The same computation applies to
$\widetilde{\widetilde{\overline{V_3'}}}$ and we get

$${\rm
rk}L_1H_2(\widetilde{\widetilde{\overline{V_3'}}})\leq 3+ 3+1=7.
$$

From this together with $(*)$, we have

$$
\begin{array}{ccc}
{\rm rk}L_1H_3(\widetilde{\widetilde{\overline{V_3}}})&\leq& {\rm
rk}H_3(W_3)+{\rm rk} L_1H_2(X)-{\rm rk} H_2(W_3)+7\\
&=& {\rm rk} H_3(W_3)-{\rm rk} H_2(W_3)+8 \quad({\rm since} \quad
L_1H_2(X)\cong \Z)
\end{array}
$$

On the other hand, we have

$$
\begin{array}{ccc}
{\rm rk}L_1H_3(\widetilde{\widetilde{\overline{V_3'}}})&\geq& {\rm
rk}H_3(W_3')+{\rm rk} L_1H_2(X')-{\rm rk} H_2(W_3')\\
&=& {\rm
rk}H_3(W_3)+{\rm rk} L_1H_2(X')-{\rm rk} H_2(W_3)\\
&=& {\rm rk} H_3(W_3)-{\rm rk} H_2(W_3)+20 \quad ({\rm since} \quad
L_1H_2(X')\cong \Z^{20})
\end{array}
$$

This shows that $L_1H_3(\widetilde{\widetilde{\overline{V_3}}})$ is
not isomorphic to $L_1H_3(\widetilde{\widetilde{\overline{V_3'}}})$.

\qe

%section 3.2**********************************************************************************************
 \ss {Application to the Case $n=3$}
{\hskip .2 in}
 With this construction, if we choose $n=3$ and $X\subset {\P^4}$
to be a general hypersurface of degree $d=5$, then $V_3=X-X\cap
{\P^3}$ and $S:=X\cap {\P^3}$ is a smooth surface in ${\P^3}$.

\medskip
\textbf{ The proof of Theorem 1.1:}
 By applying Theorem 2.1 to the pair $(X, S)$, we get
\begin{eqnarray}
 \cdots\rightarrow L_1H_3(V_3)\rightarrow L_1H_2(S)\rightarrow
L_1H_2(X)\rightarrow L_1H_{2}(V_3)\rightarrow 0.
\end{eqnarray}

The above long exact sequence $(2)$ remains exact after being
tensored with ${\mathbb{Q}}$. Note that $L_1H_2(X)\otimes
{\mathbb{Q}}\supset {\rm Griff}_1(X)\otimes {\mathbb{Q}}$ is an
infinite dimensional ${\mathbb{Q}}$-vector space by [C].  Recall
that $L_1H_2(S)$ is finitely generated since ${\dim}S=2$ (cf. [F]).
Hence $L_1H_2(V_3)\otimes {\mathbb{Q}}$ is an infinite dimensional
${\mathbb{Q}}$-vector space. By (2.2), we have $L_1H_3(V_4)\otimes
{\mathbb{Q}}\cong L_1H_2(V_3)\otimes{\mathbb{Q}}$ is an infinite
dimensional ${\mathbb{Q}}$-vector space.

\medskip
 Note that $Z_3=\overline{V_4}-{V_4}$ is defined by $(x_5\cdot
f(0,x_1,...,x_4)=0,x_0=0)$ in $\P^5$. Let
 $S'=(x_5=0)\cap(f(0,x_1,...,x_4)=0)$ in the hyperplane
 $(x_0=0)\subset {\P^5}$. It is easy to see that $S'\cong S$. Then
 $Z_3={\P^3}\cup \Sigma_p(S)$, where $\Sigma_p(S)$ means the
 joint of $S$ and the point $p=[1:0:\cdots:0]$. By applying
 Theorem 2.1 to the pair $(Z_3, \Sigma S)$, we get
\begin{eqnarray}
 \cdots\rightarrow L_1H_3(Z_3-\Sigma S)\rightarrow L_1H_2(\Sigma
S)\rightarrow L_1H_2(Z_3)\rightarrow L_1H_{2}(Z_3-\Sigma
S)\rightarrow 0.
\end{eqnarray}

Note that $Z_3-\Sigma S\cong \P^3-S$. Therefore $L_1H_2(Z_3)\otimes
{\mathbb{Q}}$ is of finite dimensional since both $L_1H_2(\sum
S)\otimes {\mathbb{Q}}\cong L_0H_0(S,{\mathbb{Q}})\cong
{\mathbb{Q}}$ ([L1]) and $L_1H_2(\P^3-S)\otimes {\mathbb{Q}}$ are.
By the same type argument, we have $L_1H_3(Z_3)\otimes {\mathbb{Q}}$
is of finite dimensional since both $L_1H_3(\sum S)\otimes {
\mathbb{Q}}\cong L_0H_1(S,{\mathbb{Q}})=0$(note that $S$ is simply
connected) and $L_1H_3(\P^3-S)\otimes {\mathbb{Q}}$ are.

\medskip
By applying Theorem 2.1 to the pair $(\overline{V_4},Z_3)$, we
have the following long exact sequence

\begin{eqnarray}
 \cdots \rightarrow L_1H_{3}({Z_3})\rightarrow
 L_1H_{3}(\overline{V_4})\rightarrow L_1H_{3}(V_{4})\rightarrow
 L_1H_{2}({Z_3})\rightarrow \cdots
\end{eqnarray}

From $(4)$, the infinite dimensionality of $L_1H_{3}({V_3})\otimes
{\mathbb{Q}}$, the finite dimensionality of $L_1H_{2}({Z_3})\otimes
{\mathbb{Q}}$ and $L_1H_{3}({Z_3})\otimes {\mathbb{Q}}$, we obtain
that $L_1H_{3}(\overline{V_4})\otimes {\mathbb{Q}}$ is  an
infinitely dimensional ${\mathbb{Q}}$-vector space. This completes
the proof of Theorem 1.1.

\qe

%
%$\bf Construction$  Let $X=(f(x_0,\cdots,x_{n+1})=0)$ be a general
%hypersurface in $\P^{n+1}$ with degree $d$, and let
%$V_{n}:=X-X\cap\{\P^{n}=(x_0=0)\}$ be the affine part, i.e.
%$V_{n}\subset \C^{n+1}$. Define $V_{n+1}:={\C^{n+1}}-V_{n}$, then
%$V_{n+1}$ can be viewed as an affine variety in ${\C^{n+2}}$
%defined by $x_{n+2}\cdot f(1,x_1,\cdots,x_{n+1})-1=0$, where
%$V_{n}=(f(1,x_1,\cdots,x_{n+1})=0)$. Now we have a sequence of
%affine spaces. Denoted by $\overline{V_{n+1}}$ the projective
%closure of $V_{n+1}$ in $\P^{n+2}$ and set $Z_{n} =
%\overline{V_{n+1}}-V_{n+1}$.

We can continue the procedure.  Set $V_{5}:={\mathbb{C}^{5}}-V_{5}$,
then $V_{5}$ can be viewed as an affine variety in
${\mathbb{C}^{6}}$ defined by $x_6\cdot(x_{5}\cdot
f(1,x_1,\cdots,x_{4})-1)-1=0$. Set $Z_{4} = \overline{V_{5}}-V_{5}$,
and so on. It can be shown in the same way that $L_1H_{3}({Z_4})$ is
finitely generated by using Theorem 2.1 and Lawson's Complex
Suspension Theorem. Note that $L_1H_{4}(V_{5})\cong L_1H_3(V_4)$ is
infinitely generated by 2.2.

By applying Theorem 2.1 to the pair $(\overline{V_5},Z_4)$, we get
the long exact sequence
$$ \rightarrow L_1H_{4}({Z_4})\rightarrow
 L_1H_{4}(\overline{V_5})\rightarrow L_1H_{4}(V_{5})\rightarrow
 L_1H_{3}({Z_4})\rightarrow \cdots$$

From these we obtain that $L_1H_{4}(\overline{V_5})$ is
infinitely generated.

{ \Prop In this construction, $L_1H_{k}(\overline{V_{k+1}})$ is
{\bf{not}} finitely generated for $k\geq 3$.}

\qe

\medskip
By the Complex Suspension Theorem [L1] ,
$L_{p+1}H_{2p+k}({\Sigma}^p\overline{V_{k+1}})\cong
L_1H_{k}(\overline {V_{k+1}})$. Therefore we get:

{\Th For integers $p$ and $k$, with $k> 0, p> 0$, we can find a
\emph{rational} projective variety $Y$,  such that $L_pH_{2p+k}(Y)$
is infinitely generated. }

 \qe
{\Rem If $k=0$ and $p>0$, there also exists projective varieties $Y$
such that $L_pH_{2p}(Y)$ is infinitely generated. This follows from
the Projective Bundle Theorem [FG] and a result of Clemens [C].}

{\Rem All the $Y$ thus constructed above are singular projective
varieties. Can one find some smooth projective variety such that the
answer to the question {\bf (Q)} is positive?  Yes, we can. The
author has constructed examples of smooth projective varieties such
{\bf (Q)} is true (cf. {\rm [H2]}). }

{\Rem Note that all $\overline {V_{k+1}}$ are singular rational
projective varieties. For smooth rational projective varieties $Y$,
$L_1H_*(Y)\otimes {\mathbb{Q}}$ are finite dimensional
${\mathbb{Q}}$-vector spaces [Pe]. The author showed  $L_1H_*(Y)$
are finitely generated abelian groups [H1]. }

%***********************************************************************************
\s {Appendix} Let $f(x_0,\cdots,x_4)$ be a general homogenous
polynomial of degree 5 and $X$ be a hypersurface of degree 6 in
$\P^5$ given by $F(x_0,\cdots,x_5):=x_5f(x_0,\cdots,x_4)-x_0^6=0$.
It is  easy to see from the proof of Lemma 3.4 that the singular
points set of X is the union of a smooth 2-dimensional variety $Y$
given by $x_0=x_5=f(x_0,x_1,\cdots,x_4)=0$ and an isolated point
defined by $\{x_0=\frac{\partial f}{\partial x_0}=\frac{\partial
f}{\partial x_1}=\frac{\partial f}{\partial
x_2}=\cdots=\frac{\partial f}{\partial x_5}=0\}$ .

Let $\sigma:\widetilde{\P^5_Y}\rightarrow \P^5$ be the blow up of
$\P^5$ along the surface $Y$ and $\tilde{X}_Y$ be the proper
transform in the blow up $\widetilde{\P^5_Y}$. Denoted by
$E=\P(N_{Y/{\P^5}})$ the exceptional divisor of the blow-up. Then
$D=E\cap \tilde{X}_Y\subset \P(N_{Y/{\P^5}})$ corresponds to the
image of the tangent cones $T_pX\subset T_p(\P^5)$ in
$\P(N_{Y/{\P^5}})$ at points $p\in Y$.

  Now $$T_pX=\bigg\{\sum_{i_0+\cdots i_5=2}\frac{\partial^2F}{\partial
  x_0^{i_0}\cdots\partial x_5^{i_5}}x_0^{i_0}\cdots x_5^{i_5}=0\bigg\}$$
  a degree 2 polynomial in $\P^5$. Directly computation shows that
  $$T_pX=\bigg\{\frac{\partial f}{\partial x_0}(p)x_0x_5+\cdots+\frac{\partial f}{\partial x_4}(p)x_4x_5)=0\bigg\}$$
  $$=(x_5=0)\cup\bigg\{\frac{\partial f}{\partial x_0}(p)x_0+\cdots+\frac{\partial f}{\partial x_4}(p)x_4)=0\bigg\}.$$
Hence $D=\tilde{X}_Y\cap E$ is a fiber bundle over $Y$ with singular
conics as fibers. Clearly, $\tilde{X}_Y$ is smooth away from $D$.
Since $D\subset E$ is a 3-dimensional variety with singular points
set $S\cong Y$, we can show that it is the only singularity on
$\tilde{X}_Y$:

{\Prop The proper transform $\tilde{X}_Y$ is a 4-dimensional variety
in $\P(N_{Y/{\P^5}})$ with singularity $S\cong Y\cup \{q\}$, where
$q$ is an isolated singular point.}

\medskip
\bp From the proof of Lemma 3.4, we see that the singular points $S$
of $X$ consist of two components. One is a smooth surface and the
other is an isolated point $q$.

Since $f$ is nonsingular on $Y=\{x_0=f(0,x_1,\cdots,x_4)=0\}$, we
have $df\neq 0$ on $Y$. Let us restrict ourselves to a neighborhood
of a point $p$ in $Y$. There, we can take the neighborhood of $p$ as
the affine space $\mathbb{C}^5$ with $p$ the origin. Hence we can
choose $y=f$ as a coordinate in the neighborhood of each point on
$Y$ since it is smooth. Locally, $Y$ is defined by $x_0=x_5=y=0$ in
$\mathbb{C}^5$. We denote it by $Y_0$. For convenience, we denote
$x_0$ by $x$, $x_5$ by $z$. The blow up
$\widetilde{(\mathbb{C}^5})_{Y_0}$ of $ \mathbb{C}^5$ along $Y_0$ is
defined by the system of equations
$$
\left\{
\begin{array}{l}
 xv=uy,\\
 xw=uz,\\
 yw=zv.
\end{array}
\right.
$$
in $\mathbb{C}^5\times P^2$, where $[u:v:w]$ is the homogenous
coordinates on $\bf P^2$. Let
$\sigma:\widetilde{(\mathbb{C}^5})_{Y_0}\rightarrow  \mathbb{C}^5$
be the map of this blowup. Then the inverse image of $X$ is given by
the following equations:
$$
\left\{
\begin{array}{l}
 x^6-yz=0,\\
 xv=uy,\\
 xw=uz,\\
 yw=zv.
\end{array}
\right.
$$
The above equations define two divisors on $\widetilde{(
\mathbb{C}^5})_{Y_0}$. One of them is the exceptional divisor $E_0$,
the intersection of $E$ with $\widetilde{(\mathbb{C}^5})_{Y_0}$ and
the other is exactly the proper transform $\widetilde{X}_{Y_0}$ of
$X_0$ in $\widetilde{(\mathbb{C}^5})_{Y_0}$, where $X_0$ is the part
of $X$ in $\mathbb{C}^5$.

 We want to show that $\widetilde{X}_{Y_0}$ is smooth away from
 $Y_0$. Now it is clear. The blow up  $\widetilde{(
 \mathbb{C}^5})_{Y_0}$ is  covered by 3 open charts: $(u\neq 0)$, $(v\neq
 0)$ and $(w\neq 0)$.

On the chart $(u\neq 0)$, we can set  $u=1$. The equations for the
inverse image of $X_0$ under $\sigma$ are given by
$$
\left\{
\begin{array}{l}
 x^6-yz=0,\\
 xv=y,\\
 xw=z,\\
 yw=zv.
\end{array}
\right.
$$
The equations $xv=y$  and  $xw=z$ imply $yw=zv$. Replacing $y$ and
$z$ by $xv$ and $xw$, respectively, we can factor $x^2$ in the first
equation $x^6-(xv)(xw)=0$. Hence the proper transform
$\widetilde{X}_{Y_0}$ are given by
$$
\left\{
\begin{array}{l}
 x^4-vw=0,\\
 xv=y,\\
 xw=z.
\end{array}
\right.
$$
and the exceptional divisor is given by

$$
\left\{
\begin{array}{l}
 x^2=0,\\
 xv=y,\\
 xw=z.
\end{array}
\right.
$$
i.e., $x=y=z=v=w=0$, which is isomorphic to $Y_0$.

It is easy to show that on the charts $(v\neq 0)$ and $(w\neq 0)$,
the proper transform $\widetilde{X}_{Y_0}$ is smooth everywhere.
This completes the proof of the proposition.

 \qe

{\Rem In fact, the 2-dimensional singularity of $X$ is of the
$A_n$-type. It can be resolved by blowing up one more time. The
isolated singularity $q$ can be resolved by one blowup.}

\begin{center}{\bf Acknowledge}\end {center}

I would like to express my gratitude to my advisor, Blaine Lawson,
for proposing this problem and for all his help.

\begin{center}{\bf Reference}\end {center}

%[AKMW] Abramovich, Dan; Karu, Kalle; Matsuki, Kenji; W\l odarczyk,
%Jaros\l aw, {\sl Torification and factorization of birational
%maps.}  J. Amer. Math. Soc. 15 (2002), no. 3, 531--572
%(electronic).

%[BS] Bloch, S.; Srinivas, V., {\sl Remarks on correspondences and
%algebraic cycles},  Amer. J. Math. 105 (1983), no. 5, 1235--1253.

\begin{enumerate}
\item[{[C]}] Clemens, H., {\sl Homological equivalence, modulo
algebraic equivalence, is not finitely generated}, I.H.E.S. Publ.
Math. 58 (1983), 19-38.

\item [{[DT]}]
Dold, A. and Thom, R.,{\sl Quasifaserungen und unendliche
symmetrische Produkte.} (German) Ann. of Math. (2) 67 1958 239--281.

\item[{[F]}] Friedlander, E. {\sl Algebraic cycles, Chow
varieties, and Lawson homology}, Comp. Math. 77(1991), 55-93.

\item[{[FG]}] Friedlander, Eric M.; Gabber, Ofer, {\sl Cycle
spaces and intersection theory}, Topological methods in modern
mathematics (Stony Brook, NY, 1991), 325--370, Publish or Perish,
Houston, TX, 1993.

\item[{[H1]}] Hu, W., {\sl Some birational invariants defined by
Lawson homology. } arXiv:math.AG/0511722.

\item[{[H2]}]Hu, W., {\sl Generalized Abel-Jacobi map on Lawson
homology.} Preprint.

%[La] Laterveer, R.,{\sl Algebraic varieties with small Chow
%groups.}  J. Math. Kyoto Univ. 38 (1998), no. 4, 673--694.

\item[{[L1]}] Lawson, H.B. Jr., {\sl Algebraic cycles and homotopy
theory}, Ann. of Math. {\bf 129}(1989), 253-291.

\item[{[L2]}] Lawson, H.B. Jr., {\sl Spaces of algebraic cycles},
pp. 137-213 in Surveys in Differential Geometry, 1995 vol.2,
International Press, 1995.

\item[{[Li]}] Lima-Filho, P., {\sl Lawson homology for
quasiprojective varieties}, Comp. Math. 84 (1992), no. 1, 1--23.

%\item[{[Pa]}] Paranjape, K.,{\sl Cohomological and cycle-theoretic
%connectivity.} Ann. Math., II. Ser. 139 (1994), 641--660

\item[{[Pe]}] Peters, C. {\sl Lawson homology for varieties with
small Chow groups and the induced filtration on the Griffiths
groups}, Math. Z. 234 (2000), no. 2, 209--223.

\item [{[V]}] Voisin, C., {\sl Th{\'e}orie de Hodge et
g{\'e}om{\'e}trie alg{\'e}brique complexe},  Cours
Sp{\'e}cialis{\'e}s, 10. Soci{\'e}t{\'e} Math{\'e}matique de
France, Paris, 2002. viii+595 pp. ISBN 2-85629-129-5

\end{enumerate}

\medskip

\noindent
Department of Mathematics,\\
Stony Brook University, SUNY,\\
Stony Brook, NY 11794-3651\\
Email:wenchuan@math.sunysb.edu

\end{document}